# Stochastic differential equations driven by fractional Brownian motions

YU-JUAN JIEN[1] and JIN MA[2]

[1]*Department of Mathematics, Purdue University, West Lafayette, IN 47907-1395, USA.
E-mail: yjien@math.purdue.edu*

[2]*Department of Mathematics, University of Southern California, Los Angeles, CA 90089, USA
and Purdue University, West Lafayette, IN 47907-1395, USA. E-mail: jinma@usc.edu*

In this paper, we study the existence and uniqueness of a class of stochastic differential equations driven by fractional Brownian motions with arbitrary Hurst parameter $H \in (0,1)$. In particular, the stochastic integrals appearing in the equations are defined in the Skorokhod sense on fractional Wiener spaces, and the coefficients are allowed to be random and even anticipating. The main technique used in this work is an adaptation of the anticipating Girsanov transformation of Buckdahn [*Mem. Amer. Math. Soc.* **111** (1994)] for the Brownian motion case. By extending a fundamental theorem of Kusuoka [*J. Fac. Sci. Univ. Tokyo Sect. IA Math.* **29** (1982) 567–597] using fractional calculus, we are able to prove that the anticipating Girsanov transformation holds for the fractional Brownian motion case as well. We then use this result to prove the well-posedness of the SDE.

## 1. Introduction

In this paper, we study the well-posedness of a class of stochastic differential equations driven by fractional Brownian motions (fBM for short) with arbitrary Hurst parameter $H \in (0,1)$ and with random coefficients that are possibly anticipating. To be more precise, we consider the following SDE:

$$X_t = X_0 + \int_0^t \sigma(s, X_s) \, \mathrm{d}B_s^H + \int_0^t b(s, X_s) \, \mathrm{d}s, \qquad (1.1)$$

where $B^H$ is a 1-dimensional fBM with parameter $H$ and $b$, $\sigma$ are measurable random fields with appropriate dimensions. At this point, we do not assume that $b$ and $\sigma$ are progressively measurable.

SDEs of this kind have been studied by many authors, mostly in the case where coefficients are deterministic, or linear (that is, $b(t,x) = b(t)x$ and $\sigma(t,x) = \sigma(t)x$, where $b(\cdot)$ and $\sigma(\cdot)$ are deterministic functions). The main difficulty is due to the fact that an fBM is neither a Markov process nor a semimartingale, except for $H = \frac{1}{2}$ (in which case $B^H$ becomes a standard Brownian motion), thus the usual stochastic calculus does not apply. As a consequence, the study of the SDE depends largely on the definitions of the stochastic integrals involved and the results vary.

We note that if $\sigma(t,x) = \sigma$ is a constant, then the SDE is of the so-called additive noise type, and the SDE involves only the Wiener integrals. In this case the path regularity of the solution does not affect the solvability directly, and the SDE can be treated as an ODE with random input. We refer to, for example, [13, 15, 18] for such case. The case where $\sigma$ is not a constant, however, is much more complicated, since the path regularity of the fBM varies with the Hurst parameter $H$ and the requirement for the path regularity of the solution varies accordingly. In particular, if $H > \frac{1}{2}$, then the paths of $B^H$ are essentially $\beta$-Hölder continuous for all $\beta < H$, hence a pathwise stochastic integral approach is quite effective (see, for example, [5, 12, 16, 17], to mention just a few).

In the general case, especially when $H < \frac{1}{2}$, the path of fBM becomes rather "rough" and the pathwise approach for stochastic integrals and the SDE becomes more difficult, therefore other definitions of stochastic integrals have been introduced. Most notable is the divergence-type integration (or Skorohod integral), which is based on the idea of Malliavin calculus for Brownian motion cases. We note that these two definitions are essentially equivalent and exchangeable (see, for example, [1, 4, 6, 7] and references cited therein). However, similar to the Brownian case, one of the main difficulties for the Skorokhod-type SDEs is that the traditional Picard iteration is no longer effective and consequently the problem becomes rather subtle when the coefficients are nonlinear and/or random. Several extended Skorokhod integrals have been defined to circumvent such difficulties, with which some special forms of SDEs have been studied (see, for example, [10, 14, 19]). However, in most of the existing literature, the diffusion coefficient $\sigma$ has to be very carefully specified so that the subtle restrictions on the stochastic integrals are satisfied. For example, it is usually assumed that $\sigma = \sigma(t,x)$ is deterministic or, even more explicitly, a linear function. In fact, to the best of our knowledge, there has not been any study of the case where both $b$ and $\sigma$ are allowed to be random and anticipating, and, at the same time, the Hurst parameter is allowed to be arbitrary.

More specifically, let us consider the following form of the SDE (1.1):

$$X_t = X_0 + \int_0^t \sigma_s X_s \, dB_s^H + \int_0^t b(s, X_s) \, ds, \qquad t \in [0,1]. \tag{1.2}$$

In the above, the stochastic integral is defined in the Skorokhod sense, $X_0$ is any $L^p$-random variable and the coefficients $\sigma$ and $b$ can be random. Our main idea is to establish a generalized version of the anticipating Girsanov theorem in the fBM setting and then to follow a scheme developed by Buckdahn [3] to attack the well-posedness of (1.2).

A major component in this method is the generalization of a fundamental theorem by Kusuoka [9] on anticipating Girsanov transformations. To be more precise, we study the following transformations $\{T_t^H, t \in [0,1]\}$ on fractional Wiener space $W$:

$$(T_t^H \omega)_\cdot = \omega_\cdot + \int_0^t K^H(\cdot, s) \sigma_s(T_s^H \omega) \, ds, \qquad \omega \in W, t \in [0,1], \tag{1.3}$$

where $K^H$ is the so-called reproducing kernel of the fBM $B^H$. We prove that, with the right choice of underlying canonical space, the probability measure induced by such a transformation is equivalent to the original one. Furthermore, similar to the Brownian case, one can also explicitly identify the Radon–Nikodym derivative of the two equivalent probability measures. Consequently, one can solve the original SDE by solving a much simpler one on a new probability space. We should note that it is the fundamental nature of this method that restricts the diffusion coefficient to being linear. However, such a restriction notwithstanding, the novelty of our result lies in the fact that the diffusion coefficient can now be random and anticipating, and the drift coefficient can even be nonlinear, which is an improvement, even compared to the original result of Buckdahn [3] in the Brownian case.

The rest of the paper is organized as follows. In Section 2, we briefly revisit some basic facts regarding fractional Brownian motion, fractional Wiener space and the Skorokhod calculus with respect to fBM. In Section 3, we study absolutely continuous transformations on fractional Wiener space and present some of their properties, and in Section 4, we revisit the Girsanov theorem of Kusuoka [9] and derive a variation of the theorem, as well as some related results. In Section 5, we present the main result on anticipating Girsanov transformation (1.3) for fBM and, finally, in Section 6, we apply these results to stochastic differential equation (1.2) and prove the existence and uniqueness of the solution.

## 2. Preliminaries

Throughout this paper, we assume that $(\Omega, \mathcal{F}, P)$ is a complete probability space and that for any $H \in (0,1)$, there exists an fBM $\{B_t^H; t \geq 0\}$, that is, a centered Gaussian process with covariance function:

$$R^H(s,t) \stackrel{\triangle}{=} E(B_s^H B_t^H) = \tfrac{1}{2}\{|s|^{2H} + |t|^{2H} - |s-t|^{2H}\}, \qquad s, t \geq 0. \tag{2.1}$$

In this paper, we assume that all processes are defined on a finite duration $[0, T]$ and, without loss of generality, we assume that $T = 1$. We shall define $I = [0, 1]$ for simplicity. Let $W \stackrel{\triangle}{=} \mathcal{C}_0(I; \mathbb{R})$ be the Banach space of continuous functions defined on $I$, null at $t = 0$ and equipped with the sup-norm. Let $\mathcal{F} \stackrel{\triangle}{=} \mathcal{B}(W)$ be the topological $\sigma$-field on $W$ and $\mu_H$ the unique probability measure on $W$ under which the canonical process $B_t^H(\omega) \stackrel{\triangle}{=} \omega_t$, $t \in I$, is an fBM. $(W, \mathcal{F}, \mu_H)$ then form a *canonical space*.

It is well known that an fBM can be represented as a Volterra-type integral of a Brownian motion. To be more precise, if $B^H$ is an fBM with $H \in (0,1)$ on $I$, then it

holds that

$$B_t^H = \int_0^t K^H(t,s) \, dB_s^{1/2}, \qquad t \in I, \tag{2.2}$$

where $K^H(t,s)$ is a non-negative function defined on $I^2$ such that $K^H(t,s) = 0$ when $s \geq t$, and it can be written explicitly in terms of Gamma and Beta functions, as well as the so-called Gaussian hypergeometric function (see [8] for details). It is clear that $R^H(s,t) = \int_0^1 K^H(t,r) K^H(s,r) \, dr$ for $s,t \in I$.

Next, we define an operator $K^H$ on $L^2(I)$ by

$$K^H f(t) \triangleq \int_0^1 K^H(t,s) f(s) \, ds, \qquad t \in I, f \in L^2(I), \tag{2.3}$$

and denote the adjoint operator of $K^H$ by $K^{H*}$. One can show that $K^{H*} \delta_{\{t\}}(s) = K^H(t,s)$, $s \in I$, where $\delta_{\{t\}}$ is the Dirac $\delta$-function at $t \in I$.

Let $W^*$ be the topological dual space of $W$ and $\mathcal{H}_H$ the associated Cameron–Martin space of $(W, \mathcal{B}(W), \mu_H)$, that is, the unique Hilbert space which is identified with its dual and is densely and continuously embedded in $W$ such that, for any $\eta \in W^*$,

$$\int_W e^{i\langle \eta, \omega \rangle} \, d\mu_H(\omega) = e^{-1/2 |\bar{\eta}|_{\mathcal{H}_H}^2}.$$

Here, $\langle \cdot, \cdot \rangle$ is the dual product $\langle \cdot, \cdot \rangle_{W^*, W}$ and $\bar{\eta} \in \mathcal{H}_H$ is the injective image of $\eta$ on $W^*$. Let us denote by $(\cdot, \cdot)_2$ and $|\cdot|_2$ the inner product and norm of $L^2(I)$, respectively. The following relations among the spaces $\mathcal{H}_H$, $L^2(I)$ and $W^*$ are useful (see [6], Theorem 3.3):

(i) $\mathcal{H}_H = K^H(L^2(I))$. More precisely, there exists $f \in L^2(I)$ for any $\tilde{f} \in \mathcal{H}_H$ such that

$$\tilde{f}(t) = K^H f(t) = \int_0^1 K^H(t,s) f(s) \, ds. \tag{2.4}$$

(ii) The scalar product on $\mathcal{H}_H$ is given by

$$(\tilde{f}, \tilde{g})_{\mathcal{H}_H} = (K^H f, K^H g)_{\mathcal{H}_H} \triangleq (f, g)_2. \tag{2.5}$$

(iii) The injection $R^H$ from $W^*$ into $\mathcal{H}_H$ can be decomposed as

$$R^H \eta = K^H(K^{H*} \eta), \qquad \eta \in W^*. \tag{2.6}$$

Since $W^*$ is continuously and densely embedded into $\mathcal{H}_H$, we define $\omega(\tilde{h}) \triangleq \lim_n \langle l_n, \omega \rangle$, where $\{l_n\}_n \subset W^*$ converges to $\tilde{h}$ in $\mathcal{H}_H$. By a slight abuse of notation, we also denote

$$\omega(h) \triangleq \omega(K^H h), \qquad h \in L^2(I), \tag{2.7}$$

when the context is clear. In what follows, we often denote $\mathcal{H}_H$ simply by $\mathcal{H}$ for a fixed $H$.

The following facts on fractional stochastic calculus can be found in [6]. We list them only for ready reference. To begin with, let $\mathcal{X}$ be a separable Hilbert space and $\mathcal{S}(\mathcal{X})$ the class of all smooth cylindrical functions $G: W \mapsto \mathcal{X}$ of the form

$$G(\omega) = g(\langle l_1, \omega\rangle, \ldots, \langle l_n, \omega\rangle)x, \qquad \omega \in W, \tag{2.8}$$

where $n \in \mathbb{N}$, $g \in C_b^\infty(\mathbb{R}^n)$, $l_k \in W^*$ for $k = 1, \ldots, n$ and $x \in \mathcal{X}$. We denote $\mathcal{S} \triangleq \mathcal{S}(\mathbb{R})$. Clearly, for any $G \in \mathcal{S}$, we can find $n \in \mathbb{N}$ and $g \in C_b^\infty(\mathbb{R}^n)$ such that

$$G(\omega) = g(\omega_{t_1}, \ldots, \omega_{t_n}), \qquad \omega \in W, 0 < t_1 < \cdots < t_n \leq 1. \tag{2.9}$$

We now define two derivatives of $G \in \mathcal{S}(\mathcal{X})$ by

$$\begin{aligned} D^{\mathcal{H}}G(\omega) &\triangleq \sum_{i=1}^n \partial_i g(\langle l_1, \omega\rangle, \ldots, \langle l_n, \omega\rangle) R^H(l_i) \otimes x, \qquad \omega \in W; \\ DG(\omega) &\triangleq \sum_{i=1}^n \partial_i g(\langle l_1, \omega\rangle, \ldots, \langle l_n, \omega\rangle) K^{H*}(l_i) \otimes x, \qquad \omega \in W. \end{aligned} \tag{2.10}$$

Then, clearly, $D^{\mathcal{H}} G \in \mathcal{H} \otimes \mathcal{X}$, but $DG \in L^2(I) \otimes \mathcal{X}$. Consequently, the directional derivatives of $G \in \mathcal{S}(\mathcal{X})$ on $\mathcal{H}$ and $L^2(I)$ are defined by

$$D_{\tilde{h}}^{\mathcal{H}} G \triangleq (D^{\mathcal{H}} G, \tilde{h})_{\mathcal{H}}, \qquad \tilde{h} \in \mathcal{H}; \qquad D_h G \triangleq (DG, h)_2, \qquad h \in L^2(I).$$

Furthermore, from (2.6) and (2.5), we have, for $\omega \in W$,

$$D^{\mathcal{H}} G(\omega) = (K^H D)G(\omega) \quad \text{and} \quad D_{\tilde{h}}^{\mathcal{H}} G(\omega) = D_h G(\omega), \qquad \text{if } \tilde{h} = K^H h. \tag{2.11}$$

We now introduce two norms in $\mathcal{S}(\mathcal{X})$ (denoting $\|\cdot\|_2$ to be the norm of $L^2(W)$),

$$\|G\|_{1,2}^{\mathcal{H}} \triangleq (\||G|_{\mathcal{X}}\|_2^2 + \||D^{\mathcal{H}} G|_{\mathcal{H} \otimes \mathcal{X}}\|_2^2)^{1/2} \quad \text{and} \quad \|G\|_{1,2} \triangleq (\||G|_{\mathcal{X}}\|_2^2 + \||DG|_{2 \otimes \mathcal{X}}\|_2^2)^{1/2},$$

and denote the closure of $\mathcal{S}(\mathcal{X})$ with respect to $\|\cdot\|_{1,2}^{\mathcal{H}}$ (resp., $\|\cdot\|_{1,2}$) by $\mathbb{D}_{\mathcal{H}}^{1,2}(\mathcal{X})$ (resp., $\mathbb{D}^{1,2}(\mathcal{X})$). The (Sobolev) spaces $\mathbb{D}_{\mathcal{H}}^{1,2}(\mathcal{X})$ and $\mathbb{D}^{1,2}(\mathcal{X})$ are then the domains of $D^{\mathcal{H}}$ and $D$, respectively. In fact, one can check that $\mathbb{D}_{\mathcal{H}}^{1,2}(\mathcal{X}) = \mathbb{D}^{1,2}(\mathcal{X})$ from (2.11) and (2.5). Finally, we define $\mathbb{D}^{1,\infty}(\mathcal{X})$ to be the space of all $G \in \mathbb{D}^{1,2}(\mathcal{X})$ such that

$$\|G\|_{1,\infty} \triangleq \||G|_{\mathcal{X}}\|_\infty \vee \||DG|_{2 \otimes \mathcal{X}}\|_\infty < \infty.$$

The following facts about the derivative $D$ are worth noting:

(i) *Chain rule.* For any random vector $G = (G_1, \ldots, G_n), n \in \mathbb{N}$, where $\{G_i\}_{i=1}^n \subset \mathbb{D}^{1,2}$, and $g \in C_b^1(\mathbb{R}^n)$, one has $g(G) \in \mathbb{D}^{1,2}$ and

$$D_t[g(G)] = \sum_{i=1}^n \partial_i g(G) D_t G_i, \qquad t \in I. \tag{2.12}$$

(ii) ([2], Proposition 2.5) If $G \in \mathbb{D}^{1,\infty} \subset \mathbb{D}^{1,2}$, then for any $\varepsilon > 0$, there exists a sequence $\{G^n\}_n \subset \mathcal{S} \subset \mathbb{D}^{1,\infty}$ which approximates $G$ in $\mathbb{D}^{1,2}$ and which satisfies, for any $n$,

$$\|G^n\|_\infty \leq \|G\|_\infty \quad \text{and} \quad \||DG^n|_2\|_\infty \leq \varepsilon + \||DG|_2\|_\infty. \tag{2.13}$$

As in the Brownian case, the Skorokhod integral with respect to an fBM is defined as the adjoint operator of the derivative operator. Namely, the integral $\delta^{\mathcal{H}}(\tilde{u})$ (resp., $\delta(u)$) is defined as the element in $L^2(W)$ such that for any $G \in \mathcal{S}$,

$$E_{\mu_H}[G\delta^{\mathcal{H}}(\tilde{u})] = E_{\mu_H}[(D^{\mathcal{H}}G, \tilde{u})_{\mathcal{H}}], \qquad \tilde{u} \in \mathcal{S}(\mathcal{H}),$$

$$(\text{resp., } E_{\mu_H}[G\delta(u)] = E_{\mu_H}[(DG, u)_2], \qquad u \in \mathcal{S}(L^2(I))).$$

From (2.11), $\delta^{\mathcal{H}}(\tilde{u}) = \delta(u)$ if $\tilde{u} = K^H u$, hence, in what follows, we often consider only $\delta$. As usual, we denote the domain of $\delta$ by $\text{Dom}(\delta)$. Then $\text{Dom}(\delta) \subset L^2(W; L^2(I))$ and a process $u \in \text{Dom}(\delta)$ if, for any $G \in \mathcal{S}$, it holds that $|E_{\mu_H}[(DG, u)_2]| \leq c\|G\|_2$, where $c$ is a constant depending on $H$ and $u$. It can be shown that $\mathbb{D}^{1,2}(L^2(I)) \subset \text{Dom}(\delta)$. We note that the spaces $\mathbb{L}^{1,2} \triangleq L^2(I; \mathbb{D}^{1,2})$ and $\mathbb{L}^{1,\infty} \triangleq L^2(I; \mathbb{D}^{1,\infty})$ are useful. They are isomorphic to $\mathbb{D}^{1,2}(L^2(I))$ and $\mathbb{D}^{1,\infty}(L^2(I))$, respectively. Thus, $\mathbb{L}^{1,2} \subset \text{Dom}(\delta)$ and one can show that

$$\|\delta(u)\|_2^2 \leq \|u\|_{1,2}^2 = \int_0^1 \|u_t\|_{1,2}^2 \, dt, \qquad u \in \mathbb{L}^{1,2}. \tag{2.14}$$

We end this section by introducing an important dense subspace of $\mathbb{L}^{1,2}$: the space of all smooth real-valued step processes, denoted by $\mathbb{L}^{\mathcal{S}}$, whose generic element is of the form

$$u_t(\omega) = g_t(\omega_{t_1}, \ldots, \omega_{t_n}), \qquad 0 < t_1 < \cdots < t_n \leq 1, (t, \omega) \in I \times W,$$

where $g: I \times \mathbb{R}^n \mapsto \mathbb{R}$ is a bounded measurable function such that $g_t(\cdot) \in C_b^\infty(\mathbb{R}^n)$ for each $t \in I$. Similar to the space $\mathbb{D}^{1,2}$, the following counterpart of (2.13) holds (see, for example, [2], Proposition 2.6): For each $u \in \mathbb{L}^{1,\infty} \subset \mathbb{L}^{1,2}$ and any $\varepsilon > 0$, there exists a sequence $\{u^n\}_n \subset \mathbb{L}^{\mathcal{S}}$ which approximates $u$ in $\mathbb{L}^{1,2}$ and is such that, for any $n$,

$$\int_0^1 \|u_s^n\|_\infty^2 \, ds \leq \int_0^1 \|u_s\|_\infty^2 \, ds \quad \text{and} \quad \int_0^1 \||Du_s^n|_2\|_\infty^2 \, ds \leq \varepsilon + \int_0^1 \||Du_s|_2\|_\infty^2 \, ds. \tag{2.15}$$

## 3. Absolutely continuous transformations on Wiener spaces

In light of the anticipating Girsanov transformation in the Brownian case, we now introduce the notion of absolutely continuous transformations on fractional Wiener spaces, this being an important component of the fractional Girsanov transformation. The difference here is that in a fractional Wiener space, such a transformation naturally involves the reproducing kernel. We shall verify that all the desired properties in [3] still hold.

Consider the fractional Wiener space $(W, \mathcal{H}, \mu) = (C_0(I; \mathbb{R}), \mathcal{H}_H, \mu_H)$ with a fixed Hurst parameter $H \in (0, 1)$. We say that a transformation $T : W \mapsto W$ is *absolutely continuous* if the image measure $\mu \circ T^{-1}$ is absolutely continuous with respect to $\mu$. The transformation $T$ is called *invertible* if there exists a transformation $A$ such that $T(A\omega) = A(T\omega) = \omega$ for all $\omega \in W$. Central to this paper is the transformation

$$T\omega = T^H \omega \stackrel{\triangle}{=} \omega + (K^H u)(\omega) = \omega_\cdot + \int_0^\cdot K^H(\cdot, r) u_r(\omega) \, dr, \tag{3.1}$$

where $u \in L^2(W; L^2(I))$ is often called the *shift process* of transformation $T$. We first state two basic properties of the transformation $T$.

**Proposition 3.1 (Lipschitz condition).** *Let $T^1$ and $T^2$ be transformations with shift processes $u^1$ and $u^2$, respectively. Assuming that either $G \in \mathcal{S}$, or $G \in \mathbb{D}^{1,\infty}$ and $T^1, T^2$ are absolutely continuous, it holds that*

$$|G(T^1 \omega) - G(T^2 \omega)| \leq \||DG|_2\|_\infty \left( \int_0^1 |u_s^1(\omega) - u_s^2(\omega)|^2 \, ds \right)^{1/2}, \qquad \mu\text{-a.e.} \tag{3.2}$$

**Proof.** We first note that if $G \in \mathcal{S}$, then the result (3.2) can be obtained directly by using the definition of derivative $D$ and (2.9).

We thus consider the case where $G \in \mathbb{D}^{1,\infty}$ and $T^i$, $i = 1, 2$, are absolutely continuous with Radon–Nikodym derivatives (or densities) $L^i, i = 1, 2$, respectively. By virtue of (2.13), there exists a sequence $\{G^n\}_n \subset \mathcal{S}$ such that $\{G^n\}_n$ converges to $G$, $\mu$-a.e. Thus,

$$E_\mu\{|G^n(T^i) - G(T^i)|\} = E_\mu\{|G^n - G|L^i\} \to 0 \qquad \text{as } n \to \infty, i = 1, 2.$$

By choosing a subsequence if necessary, we assume that the sequence $\{G^n(T^i)\}_n$ converges to $G(T^i)$ $\mu$-a.e., $i = 1, 2$. On the other hand, since $G^n \in \mathcal{S}$ for every $n$, using (3.2) and (2.13), we see that for any $\varepsilon > 0$ and $\mu$-a.e. $\omega \in W$, it holds that

$$|G^n(T^1 \omega) - G^n(T^2 \omega)| \leq (\varepsilon + \||DG|_2\|_\infty) \left( \int_0^1 |u_s^1(\omega) - u_s^2(\omega)|^2 \, ds \right)^{1/2}.$$

It follows that for any $\varepsilon > 0$ and $\mu$-a.e. $\omega \in W$, one can choose $n$ large enough such that

$$|G(T^1 \omega) - G(T^2 \omega)|$$
$$\leq |G(T^1 \omega) - G^n(T^1 \omega)| + |G^n(T^1 \omega) - G^n(T^2 \omega)| + |G^n(T^2 \omega) - G(T^2 \omega)|$$
$$\leq (\varepsilon + \||DG|_2\|_\infty) \left( \int_0^1 |u_s^1(\omega) - u_s^2(\omega)|^2 \, ds \right)^{1/2} + \varepsilon, \qquad \mu\text{-a.e.}$$

The result follows by letting $\varepsilon \to 0$ in the above. $\square$

**Proposition 3.2 (Chain rule).** *Let $G \in \mathbb{D}^{1,\infty}$ and $T$ be a transformation with shift process $u \in \mathbb{L}^{1,\infty}$. Assume that either $G \in \mathcal{S}$ or $T$ is absolutely continuous. Then $G(T) \in \mathbb{D}^{1,\infty}$ and, for any $s \in I$,*

$$D_s[G(T\omega)] = (D_s G)(T\omega) + \int_0^1 (D_r G)(T\omega)(D_s u_r)(\omega)\, dr, \qquad \mu\text{-a.e.} \qquad (3.3)$$

**Proof.** We begin by assuming that $G \in \mathcal{S}$ is as in (2.9). Then

$$G(T\omega) = g((T\omega)_{t_1}, \ldots, (T\omega)_{t_n}) = g(G_1(\omega), \ldots, G_n(\omega)), \qquad \omega \in W,$$

where $G_i(\omega) = (T\omega)_{t_i}$, $i = 1, \ldots, n$. By (3.1), (2.10) and the property of $K^{H*}$, we have

$$\begin{aligned}
D_s G_i(\omega) &= D_s[(T\omega)_{t_i}] = D_s[\langle \delta_{\{t_i\}}, \omega \rangle + \langle \delta_{\{t_i\}}, (K^H u)(\omega) \rangle] \\
&= K^{H*}\delta_{\{t_i\}}(s) + D_s[(K^{H*}\delta_{\{t_i\}}, u(\omega))_2] \\
&= K^H(t_i, s) + \int_0^{t_i} K^H(t_i, r) D_s u_r(\omega)\, dr, \quad s \in I, \omega \in W.
\end{aligned}$$

Thus, $G_i \in \mathbb{D}^{1,\infty}$, $i = 1, \ldots, n$, since $u \in \mathbb{L}^{1,\infty}$. Now, applying the chain rule (2.12), we have

$$\begin{aligned}
D_s[G(T\omega)] &= D_s[g(G_1(\omega), \ldots, G_n(\omega))] \\
&= \sum_{i=1}^n \partial_i g(G_1(\omega), \ldots, G_n(\omega)) \left( K^H(t_i, s) + \int_0^{t_i} K^H(t_i, r) D_s u_r(\omega)\, dr \right) \\
&= (D_s G)(T\omega) + \int_0^1 (D_r G)(T\omega) D_s u_r(\omega)\, dr, \qquad s \in I, \omega \in W.
\end{aligned}$$

Hence (3.3) holds when $G \in \mathcal{S}$.

To show that $G(T) \in \mathbb{D}^{1,\infty}$, we integrate the squares of both sides of equation (3.3) and then take the $L^\infty(W)$-norm. Letting $C_u = \int_0^1 \|u_r\|_{1,\infty}^2\, dr = \|u\|_{1,\infty}^2$, we have

$$\begin{aligned}
& \left\| \int_0^1 |D_s[G(T)]|^2\, ds \right\|_\infty \\
&\leq 2 \left\| \int_0^1 |(D_s G)(T)|^2\, ds \right\|_\infty \\
&\quad + 2 \left\| \int_0^1 |(D_r G)(T)|^2\, dr \cdot \int_0^1 \int_0^1 |D_s u_r|^2\, dr\, ds \right\|_\infty \\
&\leq 2(1 + C_u) \left\| \int_0^1 |(D_s G)(T)|^2\, ds \right\|_\infty = 2(1 + C_u) \left\| \int_0^1 |D_s G|^2\, ds \right\|_\infty < \infty.
\end{aligned} \qquad (3.4)$$

Note also that since $\|G(T)\|_\infty = \|G\|_\infty < \infty$, it follows that $G(T) \in \mathbb{D}^{1,\infty}$.

Now consider the general case $G \in \mathbb{D}^{1,\infty}$, but assume that $T$ is absolutely continuous with density $L$. We choose a sequence $\{G^n\}_n \subset \mathcal{S}$ satisfying (2.13) with $\varepsilon = 1$. Since $G^n(T) \in \mathbb{D}^{1,\infty}$ for any $n$ by the previous part, using a similar argument as for (3.4) and (2.13) with $\varepsilon = 1$, we obtain that for any $n$,

$$\left\| \int_0^1 |D_s[G^n(T)]|^2 \, ds \right\|_\infty \le 2(1 + C_u)\left(1 + \left\| \int_0^1 |D_s G|^2 \, ds \right\|_\infty \right).$$

Hence $\{G^n(T)\}_n$ is bounded in $\mathbb{D}^{1,\infty}$. Next, since $G^m - G^n \in \mathcal{S}$ for any $m, n \in \mathbb{N}$, replacing $G$ by $G^m - G^n$ and taking expectation instead of $L^\infty(W)$-norm in (3.4), it follows that

$$E\left\{ \int_0^1 |D_s[G^m(T) - G^n(T)]|^2 \, ds \right\} \le 2(1 + C_u) E\left\{ \int_0^1 |D_s(G^m - G^n)(T)|^2 \, ds \right\}.$$

Therefore, recalling that $\{G^n\}_n$ converges in $\mathbb{D}^{1,2}$ and letting $m, n \to \infty$, we see that

$$\|G^m(T) - G^n(T)\|_{1,2}^2 \le E|L(G^m - G^n)|^2 + 2(1 + C_u) E\left\{ L \int_0^1 |D_s(G^m - G^n)|^2 \, ds \right\} \to 0.$$

In other words, the sequence $\{G^n(T)\}_n$ converges in $\mathbb{D}^{1,2}$ and is bounded in $\mathbb{D}^{1,\infty}$, which implies that $\{G^n(T)\}_n$ converges to $G(T) \in \mathbb{D}^{1,\infty}$. Consequently, by first setting $G^n \in \mathcal{S}$ in (3.3) and then letting $n \to \infty$, we see that (3.3) holds for $G \in \mathbb{D}^{1,\infty}$, proving the proposition. $\square$

To end this section, we present the following proposition concerning the limiting behavior of the random transformation $T$.

**Proposition 3.3.** *Let $\{T^n\}_n$ be a sequence of absolutely continuous transformations with respective shift processes $\{u^n\}_n \subset L^2(W; L^2(I))$ and densities $\{L^n\}_n$. Assume that*

(i) *$\{L^n\}_n$ are uniformly integrable;*
(ii) *$\{u^n\}_n$ converges to $u$ in $L^2(W; L^2(I))$.*

*The limiting transformation $T$ defined by $T\omega \triangleq \omega + K^H u(\omega)$ is then also absolutely continuous and its density $L$ is the limit of $\{L^n\}_n$ in $L^1(W)$. Furthermore, if $\{G^n\}_n$ is a sequence of uniformly bounded random variables which converges to $G \in L^2(W)$, then the sequence $\{G^n(T^n)\}_n$ converges to $G(T)$ in $L^2(W)$ as well.*

**Proof.** Let $\{T^n\}$ and $T$ be as defined. Applying the Cauchy–Schwarz inequality on the sup-norm of $W$ and using assumption (ii), we have, for $\mu$-a.e. $\omega$,

$$|T^n\omega - T\omega|_W \triangleq \sup_{s \in I} \left| \int_0^s K^H(s,r)(u_r^n(\omega) - u_r(\omega)) \, dr \right|$$

$$\leq \sup_{s\in I}\left[\left(\int_0^s K^H(s,r)^2\,\mathrm{d}r\right)^{1/2}\right]\cdot\left(\int_0^s |u_r^n(\omega)-u_r(\omega)|^2\,\mathrm{d}r\right)^{1/2}$$
$$\to 0, \qquad n\to\infty.$$

That is, $\{T^n\}_n$ converges to $T$. Hence the sequence of measures $\{\mu\circ(T^n)^{-1}\}_n$ converges to $\mu\circ T^{-1}$, $T$ is absolutely continuous and $L$ is the limit of $\{L^n\}_n$ under the assumption (i).

To see the second half of the proposition, note that for any $n$,

$$E|G^n(T^n)-G(T)|^2 \leq 2E|G^n(T^n)-G^n(T)|^2 + 2E|G^n(T)-G(T)|^2$$
$$\leq 2E|G^n(T^n-T)|^2 + 2E|(G^n-G)L|^2.$$

Applying the result of the first part and using the fact that $\{G^n\}_n$ is uniformly bounded, the result follows immediately. $\square$

## 4. Kusuoka's theorem revisited

In this section, we turn our attention to the Girsanov transformation on fractional Wiener spaces. In light of the Brownian case, an important tool for studying such a transformation is a general theorem by Kusuoka [9]. We shall first revisit this theorem and establish some basic characterizations of the operators involved, in the context of fractional Wiener spaces.

First, let $\mathcal{X},\mathcal{X}'$ be two separable Hilbert spaces and let $\mathcal{L}_2(\mathcal{X},\mathcal{X}')$ denote the space of all Hilbert–Schmidt operators from $\mathcal{X}$ into $\mathcal{X}'$.

**Definition 4.1.** *Let $F$ be an $\mathcal{X}$-valued function defined on $W$. $F$ is called an $\mathcal{H}$-$C^1$ map if for $\mu$-a.e. $\omega\in W$, the map $h\mapsto F(\omega+h)$ is a continuous Fréchet differentiable function on $\mathcal{H}$ and its Fréchet derivative $D^{\mathcal{H}}F(\omega+\cdot):\mathcal{H}\mapsto \mathcal{L}_2(\mathcal{H},\mathcal{X})$ is continuous.*

It is known that any $\mathcal{H}$-$C^1$ map belongs to $\mathrm{Dom}(\delta^{\mathcal{H}})$ (see the corollary to Theorem 5.2 in [9]). For a generic space $V$, let $I_V$ be the identity map on space $V$ and define the Carleman–Fredholm determinant of $I_{\mathcal{X}}+B$ for $B\in\mathcal{L}_2(\mathcal{X},\mathcal{X})$ by

$$d_c(I_{\mathcal{X}}+B) = \prod_{j=1}^{\infty}(1+\lambda_j)\mathrm{e}^{-\lambda_j}, \tag{4.1}$$

where the $\lambda_j$'s are the non-zero eigenvalues of $B$, counting multiplicities. Note that $d_c(\cdot):\mathcal{L}_2(\mathcal{X},\mathcal{X})\mapsto\mathbb{R}$ is continuous. Moreover, if $B$ is a nuclear operator, then

$$d_c(I_{\mathcal{X}}+B) = \det(I_{\mathcal{X}}+B)\exp(-\operatorname{trace} B). \tag{4.2}$$

The following result of Kusuoka [9], Theorem 6.4, is crucial.

**Theorem 4.2 (Kusuoka).** *Let $K$ be an $\mathcal{H}$-$C^1$ map from $W$ to $\mathcal{H}$. Assume that for $\mu$-a.e. $\omega \in W$, the mapping $I_W + K : W \mapsto W$ is bijective and $I_{\mathcal{H}} + D^{\mathcal{H}} K(\omega) : \mathcal{H} \mapsto \mathcal{H}$ is invertible. Then $(I_W + K)^{-1} \mu(\mathrm{d}\omega) = |d(\omega; K)| \mu(\mathrm{d}\omega)$, for $\mu$-a.e. $\omega \in W$, where*

$$d(\omega; K) = d_c(I_{\mathcal{H}} + D^{\mathcal{H}} K(\omega)) \exp\{-\delta^{\mathcal{H}} K(\omega) - \tfrac{1}{2} |K(\omega)|_{\mathcal{H}}^2\}. \tag{4.3}$$

*That is, $E_\mu[G(I_W + K)|d(\cdot; K)|] = E_\mu[G]$ for any random variable $G$ on $W$.*

From the theorem, we see that the transformation involves the Carleman–Fredholm determinant of $\mathcal{L}_2(\mathcal{H}, \mathcal{H})$, as well as the $\mathcal{H}$-norm of the map. A more convenient version, which we now present, recasts the theorem in terms of $L^2(I)$ instead of $\mathcal{H}$.

**Theorem 4.3.** *Let $u : W \mapsto L^2(I)$ be a measurable mapping and $T$ a transformation defined by*

$$T\omega = \omega + (K^H u)(\omega) = \omega_\cdot + \int_0^1 K^H(\cdot, r) u_r(\omega) \, \mathrm{d}r, \qquad \omega \in W.$$

*Assume that the following conditions hold for $\mu$-a.e. $\omega \in W$:*

(i) *$T$ is bijective.*
(ii) *There exists $Du(\omega) \in L^2(I^2)$ such that for any $h \in L^2(I)$,*

(1) *$h \mapsto Du(\omega + K^H h)$ is continuous from $L^2(I)$ into $L^2(I^2)$;*
(2) *$|u_\cdot(\omega + K^H h) - u_\cdot(\omega) - (D.u(\omega), h)_2|_2 = \mathrm{o}(|h|_2)$ as $|h|_2 \to 0$;*
(3) *the mapping $I_{L^2(I)} + Du(\omega) : h \mapsto h + (Du_\cdot(\omega), h(\cdot))_2$ is invertible.*

*The measures $\mu$ and $\mu \circ T^{-1}$ are then equivalent and $A \stackrel{\triangle}{=} T^{-1}$ has the density*

$$\frac{\mathrm{d}[\mu \circ A^{-1}]}{\mathrm{d}\mu}(\omega) = |d_c(I_{L^2(I)} + Du(\omega))| \\ \times \exp\left\{-\delta u(\omega) - \frac{1}{2} |u(\omega)|_2^2\right\}, \qquad \mu\text{-a.e.}, \omega \in W. \tag{4.4}$$

**Proof.** We shall check that $K_u \stackrel{\triangle}{=} K^H u : W \mapsto \mathcal{H}$ satisfies the hypotheses in Theorem 4.2. To see this, for any $\tilde{h} \in \mathcal{H}$, let $h \in L^2(I)$ be such that $\tilde{h} = K^H h$ and $|\tilde{h}|_{\mathcal{H}} = |K^H h|_{\mathcal{H}} = |h|_2$, by virtue of (2.4) and (2.5). Hence the mapping $\tilde{h} = K^H h \mapsto D^{\mathcal{H}} K_u(\omega + \tilde{h}) = (K^H D)(K^H u)(\omega + K^H h)$ is continuous under condition (ii)(1). Moreover, by the definition of the $\mathcal{H}$-norm, one has

$$|K_u(\omega + \tilde{h}) - K_u(\omega) - D^{\mathcal{H}} K_u(\omega)(\tilde{h})|_{\mathcal{H}} = |u(\omega + K^H h) - u(\omega) - (D.u(\omega), h(\cdot))_2|_2.$$

Therefore, $K_u = K^H u$ is an $\mathcal{H}$-$C^1$ map, thanks to assumption (ii)(2). Next, note that $T = I_W + K_u$ is bijective by assumption (i). Finally, observe that if $h' = (I_{L^2(I)} + Du(\omega))(h) = h + (Du(\omega), h)_2$ for a fixed $\omega$, then

$$(I_{\mathcal{H}} + D^{\mathcal{H}} K_u(\omega))(\tilde{h}) \stackrel{\triangle}{=} \tilde{h} + (D^{\mathcal{H}} K_u(\omega), K^H h)_{\mathcal{H}} = K^H h + (D(K^H u)(\omega), h)_2 = K^H h'.$$

Thus, assumption (ii)(3) implies that $I_{\mathcal{H}} + D^{\mathcal{H}} K_u(\omega)$ is invertible on $\mathcal{H}$.

We can now apply Theorem 4.2 to conclude that the measures $\mu$ and $\mu \circ T^{-1}$ are equivalent. In order to verify the density (4.4), we first note that the operators $D^{\mathcal{H}} K_u(\omega)$ and $Du(\omega)$ have the same eigenvalues, so $d_c(I_{L^2(I)} + Du(\omega)) = d_c(I_{\mathcal{H}} + D^{\mathcal{H}} K_u(\omega))$ by definition (4.1) of the Carleman–Fredholm determinant. It therefore follows from (4.3) that

$$|d(\omega; K_u)| = |d_c(I_{\mathcal{H}} + D^{\mathcal{H}} K_u(\omega))| \exp\{-\delta^{\mathcal{H}}(K^H u)(\omega) - \tfrac{1}{2}|K^H u(\omega)|^2_{\mathcal{H}}\},$$

proving (4.4), and hence the theorem. $\square$

Now, for a given $H > 0$ and $\sigma \in \mathbb{L}^{1,\infty}$, we consider the following family of transformations $\{T^H_t, t \in I\}$ on $W$:

$$\begin{aligned}(T^H_t \omega)_s &\triangleq \omega_s + (K^H(1_{[0,t]}(\cdot)\sigma_\cdot(T^H_\cdot \omega)))_s \\ &= \omega_s + \int_0^{t \wedge s} K^H(s,r)\sigma_r(T^H_r \omega)\,\mathrm{d}r, \qquad s \in I.\end{aligned} \quad (4.5)$$

In what follows, for notational simplicity, we often drop the index $H$ from $T^H$ and $A^H$, if there is no danger of confusion. We note that the family $\{T_t\}$ is defined via differential equation (4.5) and therefore the following well-posedness result is important.

**Proposition 4.4.** *Assume $\sigma \in \mathbb{L}^{\mathcal{S}}$. Then (4.5) defines a unique family of transformations $\{T_t\}_{t \in I}$. Moreover, $T_t$ is bijective for each $t \in I$.*

**Proof.** We assume that $\sigma \in \mathbb{L}^{\mathcal{S}}$ takes the form $\sigma_t(\omega) = f_t(\omega_{t_1}, \ldots, \omega_{t_n})$, where $0 = t_0 < t_1 < \cdots < t_n \leq 1$ is any partition of $[0,1]$ and $f: I \times \mathbb{R}^n \mapsto \mathbb{R}$ is bounded and measurable such that $f_t \in C^\infty_b(\mathbb{R}^n)$ for each $t \in I$. There then exists $C_\sigma > 0$ such that for all $t \in I$,

$$|\sigma_t(\omega)| \leq C_\sigma \quad \text{and} \quad |\sigma_t(\omega) - \sigma_t(\omega')| \leq C_\sigma |\omega - \omega'|_W, \qquad \omega, \omega' \in W. \quad (4.6)$$

Consider now the following differential equation of Volterra type:

$$\xi^t_s(\omega) = \omega_s + \int_0^{t \wedge s} K^H(s,r)\sigma_r(\xi^\cdot_r(\omega))\,\mathrm{d}r, \qquad s, t \in I. \quad (4.7)$$

We show that this equation has a unique solution and that the mapping $t \mapsto \xi^t(\omega)$ is continuous in $W$ for all $\omega \in W$. To this end, let $\omega \in W$ be given and define the Picard iteration as follows. For each $t \in I$, we define $\xi^{t,0}_s(\omega) = \omega_s$ for $s \in I$, and for $n \geq 1$, we define

$$\xi^{t,n}_s(\omega) \triangleq \omega_s + \int_0^{t \wedge s} K^H(s,r)\sigma_r(\xi^{\cdot,n-1}_r(\omega))\,\mathrm{d}r, \qquad s, t \in I. \quad (4.8)$$

It is obvious that for fixed $t \in I$, $\xi^{t,n}(\omega) \in W$ for all $n$. Moreover, for $t < t'$, one has

$$|\xi^{t',n}(\omega) - \xi^{t,n}(\omega)|_W \leq C_\sigma \sup_{s \in I} \left|\int_{t \wedge s}^{t' \wedge s} K^H(s,r)\,\mathrm{d}r\right|,$$

thanks to (4.6). Consequently, the mapping $t \mapsto \xi^{t,n}(\omega)$ is continuous in $W$. This, together with the Lipschitz condition on $\sigma$ in (4.6), implies that the mapping $t \mapsto \sigma_t(\xi^{t,n-1}(\omega))$ is also continuous and hence the iteration in (4.8) is well defined. Furthermore, applying (4.6) and the Cauchy–Schwarz inequality, one can prove by induction that

$$|\xi^{t,n}(\omega) - \xi^{t,n-1}(\omega)|_W \leq \frac{C_\sigma^n t^{n/2}}{\sqrt{n!}} \leq \frac{C_\sigma^n}{\sqrt{n!}}, \qquad n \in \mathbb{N}.$$

The existence and uniqueness of the family of pathwise $W$-valued solutions $\{T_t\omega, t \in I\} \triangleq \{\xi^t(\omega), t \in I\}$ of (4.5) for $\omega \in W$ then follow from some standard argument for ordinary differential equations.

To prove the bijectiveness of $T$, we first note that an argument similar to that above also shows the well-posedness of the family $\{A_{v,t}, 0 \leq v \leq t\}$ defined by

$$\begin{aligned}(A_{v,t}\omega)_s &= \omega_s - \int_{v \wedge s}^{t \wedge s} K^H(s,r)\sigma_r(A_{r,t}\omega)\,\mathrm{d}r \\ &= \omega_s - \int_0^{t \wedge s} K^H(s,r)\sigma_r(A_{r,t}\omega)\,\mathrm{d}r + \int_0^{v \wedge s} K^H(s,r)\sigma_r(A_{r,t}\omega)\,\mathrm{d}r, \qquad \omega \in W.\end{aligned} \quad (4.9)$$

On the other hand, by (4.5), we have, for $0 \leq v \leq t$ and $\omega \in W$,

$$\begin{aligned}(T_v(A_{0,t}\omega))_s &= (A_{0,t}\omega)_s + \int_0^{v \wedge s} K^H(s,r)\sigma_r(T_r A_{0,t}\omega)\,\mathrm{d}r \\ &= \omega_s - \int_0^{t \wedge s} K^H(s,r)\sigma_r(A_{r,t}\omega)\,\mathrm{d}r + \int_0^{v \wedge s} K^H(s,r)\sigma_r(T_r A_{0,t}\omega)\,\mathrm{d}r.\end{aligned} \quad (4.10)$$

Comparing equations (4.9) and (4.10), we see that $A_{v,t} = T_v(A_{0,t})$, $0 \leq v \leq t \leq 1$, thanks to the uniqueness of $T$ and $A$. Similarly, it can be shown that $A_{v,t}(T_t) = T_v$, $0 \leq v \leq t \leq 1$. We now define $A_t \triangleq A_{0,t}$, $t \in I$. Then $T_t(A_t\omega) = T_t(A_{0,t}\omega) = \omega$ and $A_t(T_t\omega) = A_{0,t}(T_t\omega) = \omega$ for any $\omega \in W$. To wit, $T_t$ is bijective for any $t \in I$ and $A_t$ is the corresponding inverse transformation. □

It is worth noting that if $\sigma \in \mathbb{L}^{\mathcal{S}}$, then the families $\{T_t\}_{t \in I}$, $\{A_{v,t}\}_{0 \leq v \leq t \leq 1}$ and $\{A_t\}_{t \in I}$ satisfy the following relations:

$$T_v A_t = T_v A_{0,t} = A_{v,t} \quad \text{and} \quad A_{v,t} T_t = T_v, \qquad 0 \leq v \leq t \leq 1. \quad (4.11)$$

We now show that $T$ and $A$ satisfy the rest of the conditions of Theorem 4.3.

**Proposition 4.5.** *Assume that $\sigma \in \mathbb{L}^{\mathcal{S}}$. The families $\{T_t\}$ and $\{A_t\}$, defined by (4.5) and (4.9), respectively, then satisfy condition* (ii) *of Theorem 4.3.*

**Proof.** We first rewrite (4.9) as

$$(A_{v,t}\omega)_s = \omega_s - (K^H(1_{[v,t]}(\cdot)\sigma_\cdot(A_{\cdot,t}\omega)))_s, \qquad \omega \in W, 0 \leq v \leq t \leq 1, s \in I. \quad (4.12)$$

We show that, for fixed $0 \leq v \leq t \leq 1$, the mappings $u^1(\omega) \stackrel{\triangle}{=} 1_{[0,t]}\sigma(T\omega)$ and $u^2(\omega) \stackrel{\triangle}{=} 1_{[v,t]}(\cdot)\sigma.(A_{\cdot,t}\omega)$, satisfy the respective conditions for $\omega \in W$. We shall prove this only for $u^1$ (hence $T$) since the argument is similar for $u^2$ (or $A$).

Assume that $\sigma \in \mathbb{L}^{\mathcal{S}}$ takes the form $\sigma_t(\omega) = f_t(\omega_{t_1}, \ldots, \omega_{t_n})$, $0 = t_0 < t_1 < \cdots < t_n \leq 1$, as in the proof of Proposition 4.4, with $f_t \in C_b^\infty(\mathbb{R}^n)$ for each $t \in I$. We choose a complete orthonormal basis $\{e_i\}_{i \in \mathbb{N}}$ in $L^2(I)$ such that for $i = 1, \ldots, n$,

$$e_i(s) = (t_i - t_{i-1})^{-H}[K^H(t_i, s) - K^H(t_{i-1}, s)].$$

Next, we define a function $g$ on $I \times \mathbb{R}^n$ such that for each $i = 1, \ldots, n$,

$$g_t(\ldots, x_i, \ldots) = f_t\left(\ldots, \sum_{k=1}^{i}(t_k - t_{k-1})^H x_k, \ldots\right), \qquad t \in I. \tag{4.13}$$

It is obvious that $g_t \in C_b^\infty(\mathbb{R}^n)$ for any $t \in I$ and $\partial_i g$ is bounded for any $i$. Using the notation $\omega(e_i)$ for $\omega \in W$ as in (2.7), we deduce from (4.13) that

$$\sigma_t(\omega) = f_t(\omega_{t_1}, \ldots, \omega_{t_n}) = g_t(\ldots, (t_i - t_{i-1})^{-H}(\omega_{t_i} - \omega_{t_{i-1}}), \ldots)$$
$$= g_t(\ldots, \langle (K^{H*})^{-1}e_i, \omega\rangle, \ldots) = g_t(\omega(e_1), \ldots, \omega(e_n)).$$

Substituting $T_t\omega$ for $\omega$ in the above, we have

$$\sigma_t(T_t\omega) = f_t((T_t\omega)_{t_1}, \ldots, (T_t\omega)_{t_n}) = g_t((T_t\omega)(e_1), \ldots, (T_t\omega)(e_n)). \tag{4.14}$$

By the definition of $T_t$ (4.5), we see that for any $i \in \mathbb{N}$,

$$\begin{aligned}(T_t\omega)(e_i) &= \omega(e_i) + K^H(1_{[0,t]}\sigma(T\omega))(e_i) \\ &= \omega(e_i) + \langle (K^{H*})^{-1}e_i, K^H(1_{[0,t]}\sigma(T\omega))\rangle \\ &= \omega(e_i) + (e_i, 1_{[0,t]}\sigma(T\omega))_2 \\ &= \omega(e_i) + \int_0^t e_i(s)g_s((T_s\omega)(e_1), \ldots, (T_s\omega)(e_n))\,\mathrm{d}s.\end{aligned} \tag{4.15}$$

Since $(T_t\omega)(e_i) = \omega(e_i)$ when $i > n$, the mapping $\omega \mapsto (T_t\omega)(e_i)$ belongs to $\mathcal{S}$ for any $t \in I$ and $i \in \mathbb{N}$. Therefore, $\sigma(T) = \{\sigma_s(T_s), s \in I\} \in \mathbb{L}^{\mathcal{S}}$ from (4.14) and consequently $u^1(\omega) = 1_{[0,t]}\sigma(T\omega)$ satisfies both parts (1) and (2) of condition (ii) of Theorem 4.3.

It remains to check condition (ii)(3). First, for each $h \in L^2(I)$, applying the chain rule (2.12) and taking directional derivatives on both sides of equation (4.16), we have, for fixed $t \in I$ and $i \in \mathbb{N}$,

$$D_h[(T_t\omega)(e_i)]$$
$$= (e_i, h)_2 + \int_0^t e_i(s)\sum_{k=1}^{n}\partial_k g_s((T_s\omega)(e_1), \ldots, (T_s\omega)(e_n))D_h[(T_s\omega)(e_k)]\,\mathrm{d}s. \tag{4.16}$$

Since $\sigma(T)$ is of the form (4.14), equation (4.16) can be written as

$$D_h[(T_t\omega)(e_i)] = (e_i, h + (D[1_{[0,t]}\sigma(T\omega)], h)_2)_2, \qquad h \in L^2(I), i \in \mathbb{N}. \tag{4.17}$$

Now, if $h \in L^2(I)$ is such that $h + (D[1_{[0,t]}\sigma(T\omega)], h)_2 = 0$, then $D_h[(T_t\omega)(e_i)] = 0$ for any $i$ in (4.17). Therefore, $(h, e_i)_2 = 0$ for any $i$, from (4.16), and hence $h = 0$. In other words, the mapping $I_{L^2(I)} + D[1_{[0,t]}\sigma(T\omega)] : L^2(I) \mapsto L^2(I)$ is injective, and consequently bijective, which is condition (ii)(3) of Theorem 4.3. This concludes the proof. $\square$

The following proposition will play an important role in the subsequent proofs.

**Proposition 4.6.** *Assume that $\sigma \in \mathbb{L}^{\mathcal{S}}$, and let $T$ and $A$ be families of transformations of the form (4.5) and (4.12), respectively. Then*

(i) *for $0 \leq v \leq t \leq 1$, $s \in I$ and $\mu$-a.e. $\omega \in W$,*

$$D_s[\sigma_t(T_t\omega)] = (D_s\sigma_t)(T_t\omega) + \int_0^t (D_r\sigma_t)(T_t\omega)D_s[\sigma_r(T_r\omega)] \, dr,$$

$$D_s[\sigma_v(A_{v,t}\omega)] = (D_s\sigma_v)(A_{v,t}\omega) - \int_v^t (D_r\sigma_v)(A_{v,t}\omega)D_s[\sigma_r(A_{r,t}\omega)] \, dr;$$

(ii) *for $\mu$-a.e. $\omega \in W$, the Carleman–Fredholm determinants of $I_{L^2(I)} + D[1_{[0,t]}\sigma(T\omega)]$ and $I_{L^2(I)} - D[1_{[v,t]}(\cdot)\sigma.(A_{\cdot,t}\omega)]$, $0 \leq v \leq t \leq 1$, are*

$$\begin{cases} d_c(I_{L^2(I)} + D[1_{[0,t]}\sigma(T\omega)]) \\ \quad = \exp\left\{-\int_0^t \int_0^s (D_r\sigma_s)(T_s\omega)D_s[\sigma_r(T_r\omega)] \, dr \, ds\right\}, \\ d_c(I_{L^2(I)} - D[1_{[v,t]}(\cdot)\sigma.(A_{\cdot,t}\omega)]) \\ \quad = \exp\left\{-\int_v^t \int_s^t (D_r\sigma_s)(A_{s,t}\omega)D_s[\sigma_r(A_{r,t}\omega)] \, dr \, ds\right\}. \end{cases} \tag{4.18}$$

**Proof.** (i) It again suffices to prove the result for $T$. Recall that in the proof of Proposition 4.5, we actually proved that $\{\sigma_t(T_t), t \in I\} \in \mathbb{L}^{\mathcal{S}}$, provided $\sigma \in \mathbb{L}^{\mathcal{S}}$. Hence $\sigma_t \in \mathcal{S}$ and $\sigma_t(T_t) \in \mathcal{S}$ for any fixed $t$. The conclusion (i) then follows easily from Proposition 3.2.

(ii) Let $\sigma \in \mathbb{L}^{\mathcal{S}}$ and the orthonormal basis $\{e_i\}$ of $L^2(I)$ be as in the proof of Proposition 4.5. $\sigma(T)$ is then of the form (4.14). Since $D[1_{[0,t]}\sigma(T)] \in \mathcal{L}_2(L^2(I), L^2(I))$ is a nuclear operator, using (4.2) and the notation $D_i \stackrel{\triangle}{=} D_{e_i}$, for any $N \geq n$, we have that

$$d_c(I_{L^2(I)} + D[1_{[0,t]}\sigma(T\omega)])$$
$$\stackrel{\triangle}{=} \det(I_{L^2(I)} + D[1_{[0,t]}\sigma(T\omega)]) \exp(-\operatorname{trace} D[1_{[0,t]}\sigma(T\omega)])$$
$$= \det[(e_i, e_j)_2 + (e_i, D_j[1_{[0,t]}\sigma(T\omega)])_2]_{i,j=1}^N \tag{4.19}$$

$$\times \exp\left\{-\operatorname{trace}\left[\int_0^t D_j[\sigma_s(T_s\omega)]e_i(s)\,\mathrm{d}s\right]_{i,j=1}^N\right\}$$

$$= \det[D_j[(T_t\omega)(e_i)]]_{i,j=1}^N \exp\left\{-\int_0^t D_s[\sigma_s(T_s\omega)]\,\mathrm{d}s\right\},$$

where the last equality is obtained by substituting $e_j$ for $h$ in (4.17) and using the definition of derivative.

Let us now define, for each $i,j,k = 1,\ldots,N$, $P_{ij}(t) \triangleq D_j[(T_t\omega)(e_i)]$ and

$$U_{ik}(s,\omega) \triangleq \begin{cases} e_i(s)\,\partial_k g_s((T_s\omega)(e_1),\ldots,(T_s\omega)(e_n)), & 1 \leq i \leq N, 1 \leq k \leq n, \\ 0, & 1 \leq i \leq N, n+1 \leq k \leq N, \end{cases}$$

and denote the matrices $P = [P_{ij}]_{i,j=1}^N$ and $U = [U_{ik}]_{i,k=1}^N$. Substituting $e_j$ for $h$ in (4.16), we derive that for $1 \leq i,j \leq N$,

$$P_{ij}(t) = (e_i,e_j)_2 + \int_0^t \sum_{k=1}^N e_i(s)\,\partial_k g_s((T_s\omega)(e_1),\ldots,(T_s\omega)(e_n))D_j[(T_s\omega)(e_k)]\,\mathrm{d}s, \qquad t \in I.$$

That is, $P(t) = I_N + \int_0^t [UP](s)\,\mathrm{d}s$, $t \in I$, where $I_N$ is the $N \times N$ identity matrix. Hence

$$\det[D_j[(T_t\omega)(e_i)]]_{i,j=1}^N = \det(P(t)) = \exp\left\{\int_0^t \operatorname{trace} U(s)\,\mathrm{d}s\right\}$$

$$= \exp\left\{\int_0^t \sum_{i=1}^N e_i(s)\partial_i g_s((T_s\omega)(e_1),\ldots,(T_s\omega)(e_n))\,\mathrm{d}s\right\} \quad (4.20)$$

$$= \exp\left\{\int_0^t (D_s\sigma_s)(T_s\omega)\,\mathrm{d}s\right\}.$$

Therefore, combining (4.20) and Proposition 4.6(i), the determinant (4.19) becomes

$$d_c(I_{L^2(I)} + D[1_{[0,t]}\sigma(T\omega)]) = \exp\left\{\int_0^t (D_s\sigma_s)(T_s\omega)\,\mathrm{d}s\right\}\exp\left\{-\int_0^t D_s[\sigma_s(T_s\omega)]\,\mathrm{d}s\right\}$$

$$= \exp\left\{-\int_0^t \int_0^s (D_r\sigma_s)(T_s\omega)D_s[\sigma_r(T_r\omega)]\,\mathrm{d}r\,\mathrm{d}s\right\},$$

proving (ii), and hence the lemma. □

## 5. An anticipating Girsanov theorem for fBM

We are now ready to prove the main result of this paper: the anticipating Girsanov theorem for fractional Brownian motions, which can be stated as follows.

**Theorem 5.1.** *Assume that $\sigma \in \mathbb{L}^{1,\infty}$. There then exists a unique family of absolutely continuous transformations $\{T_t, t \in I\}$ such that (4.5) holds and the process $\sigma(T) = \{\sigma_t(T_t), t \in I\}$ belongs to $\mathbb{L}^{1,\infty}$. Moreover, the transformation $T_t$ is invertible for each $t \in I$ and its inverse transformation $A_t$ has density, for $\mu$-a.e. $\omega \in W$,*

$$\mathcal{L}_t(\omega) \triangleq \frac{d[\mu \circ A_t^{-1}]}{d\mu}(\omega)$$
$$= \exp\bigg\{-\delta(1_{[0,t]}\sigma(T\omega)) - \frac{1}{2}\int_0^t |\sigma_s(T_s\omega)|^2 \, ds \qquad (5.1)$$
$$- \int_0^t \int_0^s (D_r\sigma_s)(T_s\omega)D_s[\sigma_r(T_r\omega)] \, dr \, ds\bigg\}.$$

*Remark 5.2.* (i) By applying Propositions 4.4–4.6 and Theorem 4.3, we can show that Theorem 5.1 holds for $\sigma \in \mathbb{L}^{\mathcal{S}} \subset \mathbb{L}^{1,\infty}$.

(ii) Theorem 5.1 also shows that there exists a unique family of absolutely continuous transformations $\{A_{v,t}, 0 \le v \le t \le 1\}$ satisfying (4.9) with their inverse densities

$$L_{v,t}(\omega) = \exp\bigg\{\delta(1_{[v,t]}(\cdot)\sigma.(A_{\cdot,t}\omega)) - \frac{1}{2}\int_v^t |\sigma_s(A_{s,t}\omega)|^2 \, ds$$
$$- \int_v^t \int_s^t (D_r\sigma_s)(A_{s,t}\omega)D_s[\sigma_r(A_{r,t}\omega)] \, dr \, ds\bigg\}, \qquad \mu\text{-a.e.} \qquad (5.2)$$

In other words, the density of $T_t$, the inverse of $A_t$, is

$$L_t(\omega) \triangleq \frac{d[\mu \circ T_t^{-1}]}{d\mu}(\omega) = L_{0,t}(\omega), \qquad \mu\text{-a.e.}, \omega \in W. \qquad (5.3)$$

Before we prove Theorem 5.1, let us carry out a quick analysis. First, recall from (2.15) that there is a sequence $\{\sigma^n\}_n \subset \mathbb{L}^{\mathcal{S}}$ such that $\{\sigma^n\}_n$ approximates $\sigma$ in $\mathbb{L}^{1,2}$ and $\{\|\sigma^n\|_{1,\infty}\}_n$ is bounded by $\|\sigma\|_{1,\infty} \triangleq C_\sigma$. By Remark 5.2(i), we can find a family of invertible, absolutely continuous transformations $\{T_t^n, t \in I\}_n$ satisfying

$$(T_t^n\omega)_s = \omega_s + \int_0^{t\wedge s} K^H(s,r)\sigma_r^n(T_r^n\omega) \, dr, \qquad \omega \in W, s \in I.$$

Furthermore, the transformations $\{T_t^n, t \in I\}_n$ and their inverses $\{A_t^n, t \in I\}_n$ have densities $\{L_t^n\}_n$ and $\{\mathcal{L}_t^n\}_n$, respectively. In the following discussion, we shall focus only on the particular sequences $\{\sigma^n\}, \{T^n\}, \{\mathcal{L}^n\}, \{A^n\}$ and $\{L^n\}$ for the given $\sigma \in \mathbb{L}^{1,\infty}$ and we collect some important properties of the "shift processes" $\{\sigma^n(T^n)\}$ and densities $\{L^n\}$ in the following lemma.

**Lemma 5.3.** (i) *The sequence of processes $\{\sigma^n(T^n)\}$ is bounded in $\mathbb{L}^{1,\infty}$;*

(ii) *the family of densities* $\{L_t^n = \frac{\mathrm{d}\mu\circ(T_t^n)^{-1}}{\mathrm{d}\mu}, t \in I\}$ *is uniformly integrable;*

(iii) *the sequence* $\{\sigma^n(T^n)\} = \{\sigma_t^n(T_t^n), t \in I\}$ *is convergent in* $L^2(W; L^2(I))$.

**Proof.** (i) We first verify the boundedness of the derivatives of $\{\sigma^n(T^n)\}$. To this end, we apply Proposition 4.6 to $\sigma^n$ and $T^n$ to obtain

$$D_s[\sigma_t^n(T_t^n)] = (D_s\sigma_t^n)(T_t^n) + \int_0^t (D_r\sigma_t^n)(T_t^n) D_s[\sigma_r^n(T_r^n)]\,\mathrm{d}r, \qquad \mu\text{-a.e.}, s, t \in I.$$

It then follows from the Cauchy–Schwarz inequality and Gronwall's lemma that

$$\int_0^1 \||D[\sigma_t^n(T_t^n)]|_2^2\|_\infty\,\mathrm{d}t \leq 2C_\sigma^2 \exp\left\{2\int_0^1 \||D\sigma_t^n|_2^2\|_\infty\,\mathrm{d}t\right\} \leq 2C_\sigma^2 \exp\{2C_\sigma^2\}. \quad (5.4)$$

Combining (5.4) with the fact that $\int_0^1 \|\sigma_t^n(T_t^n)\|_\infty^2\,\mathrm{d}t \leq C_\sigma^2$ for every $n$, we have

$$\int_0^1 \|\sigma_t^n(T_t^n)\|_{1,\infty}^2\,\mathrm{d}t \leq C_\sigma^2 + 2C_\sigma^2 \exp\{2C_\sigma^2\}, \qquad n \in \mathbb{N}, \quad (5.5)$$

proving (i).

To see (ii), we choose $\phi(x) = x|\ln x|$. It then suffices to show that

$$\sup_{t \in I, n \in \mathbb{N}} E\{\phi(L_t^n)\} = \sup_{t \in I, n \in \mathbb{N}} E\{L_t^n|\ln L_t^n|\} < \infty. \quad (5.6)$$

From the definition of density, we have $(L_t^n)^{-1} = \mathcal{L}_t^n(A_t^n)$, $t \in I$, for any $n \in \mathbb{N}$. Therefore,

$$E\{L_t^n|\ln L_t^n|\} = E\{L_t^n|\ln \mathcal{L}_t^n(A_t^n)|\} = E\{|\ln \mathcal{L}_t^n(A_t^n T_t^n)|\} = E\{|\ln \mathcal{L}_t^n|\}.$$

Using (5.1), we obtain

$$E\{|\ln \mathcal{L}_t^n|\} \leq E\{|\delta(1_{[0,t]}\sigma^n(T^n))|\} + E\left\{\frac{1}{2}\int_0^t |\sigma_s^n(T_s^n)|^2\,\mathrm{d}s\right\}$$
$$+ E\left\{\left|\int_0^t \int_0^s (D_r\sigma_s^n)(T_s^n) D_s[\sigma_r^n(T_r^n)]\,\mathrm{d}r\,\mathrm{d}s\right|\right\} \stackrel{\triangle}{=} I_1 + I_2 + I_3.$$

We shall find the upper bounds for $I_i$, $i = 1, 2, 3$. First, note that for each $n$, and $t \in I$,

$$I_1 \leq \|\delta(1_{[0,t]}\sigma^n(T^n))\|_2 \leq \|\sigma^n(T^n)\|_{1,2} \leq \|\sigma^n(T^n)\|_{1,\infty} \leq (C_\sigma^2 + 2C_\sigma^2 \exp\{2C_\sigma^2\})^{1/2},$$

where the second inequality is obtained by applying (2.14) and the last inequality by applying (5.5). Next, since $\||\sigma^n|_2\|_\infty \leq C_\sigma$ for any $n$, $I_2 \leq \frac{1}{2}C_\sigma^2$ for all $t$. Finally, applying the Cauchy–Schwarz inequality, we have

$$I_3 \leq \left(\int_0^1 \int_0^1 \|(D_r\sigma_s^n)(T_s^n)\|_\infty^2\,\mathrm{d}r\,\mathrm{d}s\right)^{1/2} \left(\int_0^1 \int_0^1 \|D_s[\sigma_r^n(T_r^n)]\|_\infty^2\,\mathrm{d}r\,\mathrm{d}s\right)^{1/2}$$
$$\leq C_\sigma \cdot \sqrt{2}C_\sigma \exp\{C_\sigma^2\} = \sqrt{2}C_\sigma^2 \exp\{C_\sigma^2\}, \qquad t \in I,$$

where the last inequality is due to (5.4). Consequently, $E\{|\ln \mathcal{L}_t^n|\} = E\{L_t^n|\ln L_t^n|\}$ is bounded, uniformly in $t \in I$ and $n$, and (5.6) (and hence (ii)) follows.

It remains to prove (iii). By using the Cauchy–Schwarz inequality and applying Proposition 3.1, we see that for $m, n \in \mathbb{N}$,

$$E\bigg\{\int_0^t |\sigma_s^m(T_s^m) - \sigma_s^n(T_s^n)|^2 \, \mathrm{d}s\bigg\}$$

$$\leq 2E\bigg\{\int_0^t |\sigma_s^m(T_s^m) - \sigma_s^n(T_s^m)|^2 \, \mathrm{d}s\bigg\} + 2E\bigg\{\int_0^t |\sigma_s^n(T_s^m) - \sigma_s^n(T_s^n)|^2 \, \mathrm{d}s\bigg\}$$

$$\leq 2E\bigg\{\int_0^t |\sigma_s^m - \sigma_s^n|^2 L_s^m \, \mathrm{d}s\bigg\} + 2E\bigg\{\int_0^t \bigg(\||D\sigma_s^n|_2^2\|_\infty \int_0^s |\sigma_r^m(T_r^m) - \sigma_r^n(T_r^n)|^2 \, \mathrm{d}r\bigg) \, \mathrm{d}s\bigg\}.$$

First applying Gronwall's lemma, then letting $t = 1$ (noting that $\{L^n\}$ is uniformly integrable, thanks to part (ii) above) we obtain that

$$E\bigg\{\int_0^1 |\sigma_s^m(T_s^m) - \sigma_s^n(T_s^n)|^2 \, \mathrm{d}s\bigg\}$$

$$\leq 2E\bigg\{\int_0^1 |\sigma_s^m - \sigma_s^n|^2 L_s^m \, \mathrm{d}s\bigg\} \exp\bigg\{2\int_0^1 \||D\sigma_s^n|_2^2\|_\infty \, \mathrm{d}s\bigg\}$$

$$\leq 2\exp\{2C_\sigma^2\} E\bigg\{\int_0^1 |\sigma_s^m - \sigma_s^n|^2 L_s^m \, \mathrm{d}s\bigg\} \longrightarrow 0$$

as $m, n \to \infty$ since $\{\sigma^n\}$ converges in $\mathbb{L}^{1,2}$. Thus, $\{\sigma^n(T^n)\}$ is a Cauchy sequence and it converges in $L^2(W; L^2(I))$. This completes the proof. □

**Proof of Theorem 5.1.** First, we recall from Lemma 5.3(i) and (iii) that the sequence $\{\sigma^n(T^n)\}$ is bounded in $\mathbb{L}^{1,\infty}$ and convergent in $\mathbb{L}^{1,2}$. Let $\bar{\sigma} \in \mathbb{L}^{1,\infty}$ be the limit process of shift processes $\{\sigma^n(T^n)\}$ and define the family of transformations $\{T_t, t \in I\}$ by

$$(T_t\omega)_s \stackrel{\triangle}{=} \omega_s + \int_0^{t \wedge s} K^H(s, r) \bar{\sigma}_r(\omega) \, \mathrm{d}r, \qquad \mu\text{-a.e.} \tag{5.7}$$

Proposition 3.3, together with Lemma 5.3(ii), then shows that $T_t$ is absolutely continuous and $\{\sigma_t^n(T_t^n)\}$ converges to $\sigma_t(T_t)$ in $L^2(W)$ for $t \in I$. In other words, $\bar{\sigma}_t = \sigma_t(T_t)$, $\mu$-a.e., for any $t \in I$ and $\{T_t, t \in I\}$ satisfies

$$(T_t\omega)_s = \omega_s + \int_0^{t \wedge s} K^H(s, r) \sigma_r(T_r\omega) \, \mathrm{d}r, \qquad \mu\text{-a.e.} \tag{5.8}$$

The uniqueness of $\{T_t, t \in I\}$ can be obtained easily by using Gronwall's lemma.

Next, recall from Remark 5.2(ii) that there exist transformations $\{A^n_{v,t}, 0 \leq v \leq t \leq 1\}$ satisfying

$$(A^n_{v,t}\omega)_s = \omega_s - \int_{v \wedge s}^{t \wedge s} K^H(s,r)\sigma^n_r(A^n_{r,t}\omega)\,dr, \qquad \mu\text{-a.e.}$$

for every $n$. Following an argument similar to that of Lemma 5.3(ii), one can show that the family of densities $\{\mathcal{L}^n_{v,t} \stackrel{\triangle}{=} \frac{d\mu \circ (A^n_{v,t})^{-1}}{d\mu}, 0 \leq v \leq t \leq 1\}$ is uniformly integrable. Hence, using the same argument as was used for (5.7), one can show that the limit of $\{\sigma^n_v(A^n_{v,t})\}$ in $\mathbb{L}^{1,2}$ exists for any $t \in I$ and denote this by $\{\bar{\sigma}_{v,t}, 0 \leq v \leq t\} \in \mathbb{L}^{1,\infty}$. Furthermore, similar to (5.8), one can also argue that the transformations defined by

$$(A_{v,t}\omega)_s \stackrel{\triangle}{=} \omega_s - \int_{v \wedge s}^{t \wedge s} K^H(s,r)\bar{\sigma}_{r,t}(\omega)\,dr, \qquad \mu\text{-a.e.} \tag{5.9}$$

are unique and absolutely continuous, and that $\{A_{v,t}, 0 \leq v \leq t \leq 1\}$ satisfy

$$(A_{v,t}\omega)_s = \omega_s - \int_{v \wedge s}^{t \wedge s} K^H(s,r)\sigma_r(A_{r,t}\omega)\,dr, \qquad \mu\text{-a.e.}$$

We now check that for each $t \in I$, the transformation $A_t \stackrel{\triangle}{=} A_{0,t}$ is the inverse of $T_t$. Note that the $\sigma^n(T^n)$'s are bounded in $\mathbb{L}^{1,\infty}$ by Lemma 5.3. Applying Proposition 3.3 to $\{\sigma^n_v(T^n_v)\}$ and $\{A^n_t\}$, and using (5.7) as well as (4.11), we conclude that $\{\sigma^n_v(A^n_{v,t}) = \sigma^n_v(T^n_v A^n_t), 0 \leq v \leq t\}$ converges to $\bar{\sigma}_v(A_t)$ for any $t \in I$. As a result, $\bar{\sigma}_{v,t} = \bar{si}_v(A_t)$ in $L^2(W)$, $0 \leq v \leq t \leq 1$. By setting $v = 0$ and substituting $T_t\omega$ for $\omega$ in (5.9), we have

$$(A_t(T_t\omega))_s = (T_t\omega)_s - \int_0^{t \wedge s} K^H(s,r)\bar{\sigma}_{r,t}(T_t\omega)\,dr$$

$$= (T_t\omega)_s - \int_0^{t \wedge s} K^H(s,r)\bar{\sigma}_r(A_t T_t\omega)\,dr = \omega_s, \qquad \mu\text{-a.e.}$$

Similarly, one shows that $T_t(A_t\omega) = \omega$, $\mu$-a.e. To wit, $A_t$ is the inverse of $T_t$ for any $t \in I$.

It remains to compute the densities of the transformations $\{A_t, t \in I\}$. From Remark 5.2(i), we see that for any $n$,

$$\mathcal{L}^n_t(\omega) = \exp\Bigg\{-\delta(1_{[0,t]}\sigma^n(T^n\omega)) - \frac{1}{2}\int_0^t |\sigma^n_s(T^n_s\omega)|^2\,ds$$

$$- \int_0^t \int_0^s (D_r\sigma^n_s)(T^n_s\omega)D_s[\sigma^n_r(T^n_r\omega)]\,dr\,ds\Bigg\}, \qquad \mu\text{-a.e.}$$

Since the sequence $\{\mathcal{L}^n_t, t \in I\}$ is uniformly integrable, by Proposition 3.3, it converges to the right-hand side of (5.1), which is the density of $\{A_t, t \in I\}$. This completes the proof. $\square$

## 6. Stochastic differential equations driven by fBM

In this section, we fix the Hurst parameter $H \in (0,1)$ and denote $B^H$ by $B$ for simplicity. Also, we will use the conventional notation $\int_0^t u_s \, dB_s$ to denote $\delta(1_{[0,t]} u)$ if $1_{[0,t]} u \in \text{Dom}(\delta)$ for $t \in I$. Note that since $H$ can be arbitrary, the Skorokhod integral should often be understood in the general sense defined by [11].

We consider the stochastic differential equation in the Skorokhod sense

$$X_t = X_0 + \int_0^t \sigma_s X_s \, dB_s + \int_0^t b(s, X_s) \, ds, \qquad t \in I, \tag{6.1}$$

where $X_0 \in L^p(W)$ for some $p \geq 2$, $\sigma \in \mathbb{L}^{1,\infty}$ and $b: I \times \mathbb{R} \times W \mapsto \mathbb{R}$ is a measurable function satisfying the following condition for $\mu$-a.e. $\omega \in W$:

**H1.** There exist an integrable function $\gamma_t \geq 0$ on $I$ and a constant $M > 0$ such that

(i) $\int_0^1 \gamma_t \, dt \leq M$ and $|b(t, 0, \omega)| \leq M$ for any $t \in I$;
(ii) $|b(t, x, \omega) - b(t, y, \omega)| \leq \gamma_t |x - y|$ for all $x, y \in \mathbb{R}$, $t \in I$.

Our plan of attack is as follows. Since $\sigma \in \mathbb{L}^{1,\infty}$, by Theorem 5.1, we know that the family of transformations on $W$ defined by

$$(T_t \omega)_s \stackrel{\triangle}{=} \omega_s + \int_0^{t \wedge s} K^H(s,r) \sigma_r(T_r \omega) \, dr, \qquad s, t \in I, \mu\text{-a.e.},$$

exists and is unique. Also, there is a corresponding family of inverse transformations $\{A_t, t \in I\}$. Using the transformations $T$ and $A$, we can define a change of probability measure using the density of $T$, and we show that the SDE becomes a much simpler one under the new probability measure. To carry out this scheme, we need the following lemma regarding the temporal derivatives of the image processes under transformations $T$ and $A$.

**Lemma 6.1.** (i) *Suppose that $F = \{F_t, t \in I\} \in \mathbb{L}^{\mathcal{S}}$ and the mapping $t \mapsto F_t(\cdot)$ is differentiable. Then $\{F_t(T_t), t \in I\}$ is differentiable with respect to $t$ and it holds that*

$$\frac{d}{dt}[F_t(T_t)] = \left(\frac{d}{dt} F_t\right)(T_t) + \sigma_t(T_t)(D_t F_t)(T_t), \qquad \mu\text{-a.e.} \tag{6.2}$$

(ii) *For any $G \in \mathcal{S}$, the mapping $t \mapsto G(A_t)$ is differentiable and it holds that*

$$\frac{d}{dt} G(A_t) = -\sigma_t D_t[G(A_t)], \qquad \mu\text{-a.e.} \tag{6.3}$$

**Proof.** Assume that $F \in \mathbb{L}^{\mathcal{S}}$ takes the form $F_t(\omega) = g_t(\langle \delta_{\{t_1\}}, \omega \rangle, \ldots, \langle \delta_{\{t_n\}}, \omega \rangle)$, $0 < t_1 < \cdots < t_n \leq 1$, where $g$ is a bounded measurable function on $I \times \mathbb{R}^n$ with $g_t \in C_b^\infty(\mathbb{R}^n)$ for each $t \in I$ and $\delta_{\{\cdot\}}$ is the Dirac $\delta$-function. Since $F_t(T_t\omega) = g_t(\langle \delta_{\{t_1\}}, T_t\omega \rangle, \ldots, \langle \delta_{\{t_n\}}, T_t\omega \rangle)$, we need only check the differentiability of $\langle \delta_{\{t_i\}}, T_t\omega \rangle$ for any $i$. Indeed,

$$\frac{d}{dt}\langle \delta_{\{t_i\}}, T_t\omega \rangle = \frac{d}{dt}\left[\omega_{t_i} + \int_0^{t \wedge t_i} K^H(t_i, r)\sigma_r(T_r\omega) \, dr\right]$$
$$= K^H(t_i, t)\sigma_t(T_t\omega) = K^{H*}(\delta_{\{t_i\}})(t)\sigma_t(T_t\omega), \qquad \mu\text{-a.e.}$$

Now applying the chain rule (2.12) and the definition of derivative operator (2.10), we have

$$\frac{d}{dt}[F_t(T_t\omega)] = \left(\frac{d}{dt}F_t\right)(T_t\omega) + \sum_{i=1}^n \partial_i g_t(\langle \delta_{\{t_1\}}, T_t\omega \rangle, \ldots, \langle \delta_{\{t_n\}}, T_t\omega \rangle) K^{H*}(\delta_{\{t_i\}})(t)\sigma_t(T_t\omega)$$
$$= \left(\frac{d}{dt}F_t\right)(T_t\omega) + (D_t F_t)(T_t\omega)\sigma_t(T_t\omega), \qquad \mu\text{-a.e.},$$

proving (i). To prove (ii), we first note that by following the same argument as was used for Proposition 4.5, one can show that the process $\{G(A_t), t \in I\} \in \mathbb{L}^{\mathcal{S}}$ for $G \in \mathcal{S}$. Therefore, from part (i) above, we obtain

$$0 = \frac{d}{dt}G = \frac{d}{dt}G(A_t T_t) = \left(\frac{d}{dt}G(A_t)\right)(T_t) + \sigma_t(T_t)D_t[G(A_t)](T_t), \qquad \mu\text{-a.e.}$$

Thus, part (ii), and hence the lemma, follows. $\square$

As in Theorem 5.1, we denote the densities of $A_t$ and $T_t$ by $\mathcal{L}_t$ and $L_t$, respectively. Let us now consider the following ordinary differential equation for any fixed $\omega \in W$:

$$Z_t(\omega, x) = x + \int_0^t L_s^{-1}(T_s\omega)b(s, L_s(T_s\omega)Z_s(\omega, x), T_s\omega) \, ds, \qquad x \in \mathbb{R}, t \in I. \quad (6.4)$$

It is known from ODE theory that under Assumption (H1), the unique solution $Z_t(\omega, x)$, $t \geq 0$, depends continuously on the initial state $x$. Thus, the mapping $(t, \omega) \mapsto Z_t(\omega, X_0(\omega))$ defines a measurable process. Let us now set

$$X_t = L_t Z_t(A_t, X_0(A_t)), \qquad t \in I. \quad (6.5)$$

The main result of this paper is the following theorem.

**Theorem 6.2.** *The process $\{X_t, t \in I\}$ in (6.5) satisfies $1_{[0,\tau]}\sigma X \in \text{Dom}(\delta)$ for all $\tau \in I$ and $X \in L^2(W; L^2(I))$ is the unique solution of the SDE (6.1).*

**Proof.** *Existence.* We will show that $1_{[0,\tau]}\sigma X \in \text{Dom}(\delta)$ for $\tau \in I$ and that SDE (6.1) holds. To this end, let $G \in \mathcal{S}$ and denote $Z_t(\cdot, X_0(\cdot))$ by $Z_t(X_0)$. Using (6.5), we have

$$E\left\{G \int_0^1 1_{[0,\tau]}(t)\sigma_t X_t \, dB_t\right\} = E\left\{\int_0^\tau \sigma_t X_t D_t G \, dt\right\}$$

$$= E\left\{\int_0^\tau \sigma_t L_t Z_t(A_t, X_0(A_t)) D_t G \, dt\right\} \tag{6.6}$$

$$= E\left\{\int_0^\tau \sigma_t(T_t) Z_t(X_0) D_t G(T_t) \, dt\right\}.$$

Applying Lemma 6.1(i) and integration by parts, (6.6) becomes

$$E\left\{\int_0^\tau \sigma_t(T_t) Z_t(X_0) D_t G(T_t) \, dt\right\} = E\left\{\int_0^\tau Z_t(X_0) \frac{d}{dt} G(T_t) \, dt\right\}$$

$$= E\left\{Z_\tau(X_0) G(T_\tau) - Z_0(X_0) G \right. \tag{6.7}$$

$$\left. - \int_0^\tau \left(\frac{d}{dt} Z_t(X_0)\right) G(T_t) \, dt\right\}.$$

Next, using ODE (6.4) as well as the fact that $L_t^{-1}(T_t) = \mathcal{L}_t$ is the density of $A_t$, (6.7) yields that

$$E\left\{Z_\tau(X_0) G(T_\tau) - Z_0(X_0) G - \int_0^\tau L_t^{-1}(T_t) b(t, L_t(T_t) Z_t(X_0), T_t) G(T_t) \, dt\right\}$$

$$= E\{L_\tau Z_\tau(A_\tau, X_0(A_\tau)) G\} - E\{Z_0(X_0) G\} - E\left\{\int_0^\tau b(t, L_t Z_t(A_t, X_0(A_t))) G \, dt\right\}$$

$$= E\{X_\tau G\} - E\{X_0 G\} - E\left\{\int_0^\tau b(t, X_t) G \, dt\right\} = E\left\{G\left(X_\tau - X_0 - \int_0^\tau b(t, X_t) \, dt\right)\right\}.$$

This, together with (6.6), leads to the fact that for any $G \in \mathcal{S}$,

$$E\left\{G \int_0^1 1_{[0,\tau]}(t) \sigma_t X_t \, dB_t\right\} = E\left\{G\left(X_\tau - X_0 - \int_0^\tau b(t, X_t) \, dt\right)\right\}.$$

Since $X_\tau - X_0 - \int_0^\tau b(t, X_t) \, dt$ is square-integrable, we deduce that $\{1_{[0,\tau]} \sigma X, \tau \in I\}$ belong to $\text{Dom}(\delta)$ and $X$ satisfies (6.1).

*Uniqueness.* Let $Y \in L^2(W; L^2(I))$, where $1_{[0,t]} \sigma Y \in \text{Dom}(\delta)$ for all $t \in I$, be any solution of equation (6.1), that is,

$$Y_t = X_0 + \int_0^t \sigma_s Y_s \, dB_s + \int_0^t b(s, Y_s) \, ds, \qquad t \in I. \tag{6.8}$$

We consider a fixed $t \in I$ and a random variable $G \in \mathcal{S}$. Multiplying both sides of (6.8) by $G(A_t)$ and taking expectations, it becomes

$$E\{Y_t G(A_t)\} = E\{Y_0 G(A_t)\} + E\left\{\int_0^t D_s[G(A_t)] \sigma_s Y_s \, ds\right\} + E\left\{\int_0^t b(s, Y_s) G(A_t) \, ds\right\}.$$

Since $G(A_t) = G(A_s) - \int_s^t \sigma_r D_r[G(A_r)] \, dr$ for any $s \in [0, t]$ by Lemma 6.1(ii), we obtain

$$\begin{aligned}
E\{Y_t G(A_t)\} &= E\{Y_0 G\} - E\left\{Y_0 \int_0^t \sigma_r D_r[G(A_r)] \, dr\right\} \\
&\quad + E\left\{\int_0^t D_s[G(A_s)] \sigma_s Y_s \, ds\right\} - E\left\{\int_0^t D_s\left[\int_s^t \sigma_r D_r[G(A_r)] \, dr\right] \sigma_s Y_s \, ds\right\} \\
&\quad + E\left\{\int_0^t b(s, Y_s) G(A_s) \, ds\right\} - E\left\{\int_0^t b(s, Y_s) \int_s^t \sigma_r D_r[G(A_r)] \, dr \, ds\right\} \\
&= E\{Y_0 G\} + E\left\{\int_0^t D_s[G(A_s)] \sigma_s Y_s \, ds\right\} + E\left\{\int_0^t b(s, Y_s) G(A_s) \, ds\right\} \\
&\quad - E\left\{\int_0^t \sigma_r D_r[G(A_r)] Y_0 \, dr\right\} - E\left\{\int_0^t \int_0^r D_s[\sigma_r D_r[G(A_r)]] \sigma_s Y_s \, ds \, dr\right\} \\
&\quad - E\left\{\int_0^t \sigma_r D_r[G(A_r)] \int_0^r b(s, Y_s) \, ds \, dr\right\}.
\end{aligned} \quad (6.9)$$

Here, the last equality is due to Fubini's theorem. Now, by definition of the Skorokhod integral,

$$E\left\{\int_0^t \int_0^r D_s[\sigma_r D_r[G(A_r)]] \sigma_s Y_s \, ds \, dr\right\} = E\left\{\int_0^t \sigma_r D_r[G(A_r)] \int_0^r \sigma_s Y_s \, dB_s \, dr\right\}.$$

Note that because the density of $A_t$ is $\mathcal{L}_t = L_t^{-1}(T_t)$ and $Y$ satisfies (6.8), (6.9) can be rewritten as

$$\begin{aligned}
E\{L_t^{-1}(T_t) Y_t(T_t) G\} &= E\{Y_0 G\} + E\left\{\int_0^t D_s[G(A_s)] \sigma_s Y_s \, ds\right\} \\
&\quad + E\left\{\int_0^t b(s, Y_s) G(A_s) \, ds\right\} - E\left\{\int_0^t \sigma_r D_r[G(A_r)] Y_r \, dr\right\} \\
&= E\{Y_0 G\} + E\left\{\int_0^t b(s, Y_s) G(A_s) \, ds\right\} \\
&= E\{Y_0 G\} + E\left\{\int_0^t L_s^{-1}(T_s) b(s, Y_s(T_s)) G \, ds\right\}.
\end{aligned}$$

Since the smooth random variable $G$ is arbitrary, we have

$$\begin{aligned}
L_t^{-1}(T_t) Y_t(T_t) &= Y_0 + \int_0^t L_s^{-1}(T_s) b(s, Y_s(T_s)) \, ds \\
&= Y_0 + \int_0^t L_s^{-1}(T_s) b(s, L_s(T_s) L_s^{-1}(T_s) Y_s(T_s)) \, ds, \qquad \mu\text{-a.e.}
\end{aligned}$$

That is, $L_t^{-1}(T_t)Y_t(T_t)$ is a solution of equation (6.4). By the uniqueness of the ODE, we must have $L_t^{-1}(T_t)Y_t(T_t) = Z_t(\cdot, Y_t)$. Consequently,

$$Y_t = L_t Z_t(A_t, Y_0(A_t)) = X_t, \qquad \mu\text{-a.e.},$$

which is the unique solution of SDE (6.1). This completes the proof. □

## Acknowledgments

The second author was supported in part by NSF Grants #0505427 and #0835051.